\documentclass[a4paper,12pt]{article}

\usepackage{amsmath}
\usepackage{amsthm}
\usepackage{amsfonts}

\parskip\smallskipamount
\newtheorem*{theorem*}{Theorem}

\newtheorem{thm}{Theorem}
\newtheorem{lem}[thm]{Lemma}

\newtheorem{prop}[thm]{Proposition}

\newcommand{\pr}{\mathbb{P}}

\newcommand{\Z}{\mathbb{Z}}
\newcommand{\M}{\mathcal{M}}
\newcommand{\Zm}{\Z^{\mbox{\bf{-}}}}
\newcommand{\Zmg}{Z^{\mbox{\bf{-}}}}

\title{Regeneration of extremal particles for one-dimensional contact processes}
\author{A. Tzioufas \vspace{-6mm} \footnote{{\textit{Departamento de Matematica, 
Universidad de Buenos Aires, 
Ciudad Universitaria, Capital Federal 
C1428EGA - Argentina}
}}}
\date{}

\begin{document}
\maketitle 
\vspace{-10mm}
\begin{abstract}
A new, conceptual proof approach for establishing the existence of regenerative space-time points for symmetric, translation invariant, finite-range interaction contact processes on survival is shown.  The proof is elementary, complements the original one, and employs symmetry-based coupling arguments and a new consequence of convergence to equilibrium of the process in order to circumvent the original block construction.

\end{abstract}
\vspace{2mm}

The extremal (rightmost and leftmost) particles of supercritical nearest and non-nearest (translation invariant and finite range) neighbour contact processes are known to obey a central limit theorem; cf.\ \cite{GP87, K89} and \cite{MS00} respectively. In the former case, relying on the consequence of the nearest-neighbours interaction assumption that the rightmost descendant of a surviving particle (i.e.\ has infinite descendancy) cannot be surpassed by descendants of particles initially to its left, the conjectured result was established in \cite{GP87} and later proved again in a simple way in \cite{K89} by means of a seminal restart coupling argument that considers random space-time points  at which the rightmost particle is conditioned to survive -- these points occur with positive probability as the process is supercritical and, owing to the consequence of the basic coupling mentioned, possess the desired regenerative property. As noted in \cite{MS00} the key to the extension of this argument to the non-nearest neighbors case is showing that there still is a positive probability for a particle which has infinite descendancy not to be surpassed for all times. 


A new approach for proving this by exploiting the symmetric interaction assumption is shown  (see Theorem \ref{Rr} below for precise formulation). This approach provides with useful insights by relying on a new simple consequence of a mode of convergence to equilibrium for the process (v.i.\ Theorem  \ref{gen}) applied for the process constrained on the half line (v.i.\ Proposition \ref{cons}), and then evoking equality of critical values for survival on the half and the whole line. In addition, it is much simpler than the known one in the sense that it does not require any renormalization group type arguments, besides that for the proofs of the before-mentioned invoked results. We finally show that this result suffices for reproving the existence of random regeneration times, i.e. times at which the process forgets its past, for the extremal particles (v.i. Theorem 1) by making use of the previously mentioned observation of \cite{MS00}. We thus show here the i.i.d.\ nature of the growth of the diameter of the process on survival, and complement the original proof (see Remark in $\mathsection$\textit{2.}) We also note that the construction in \cite{MS00} allows for obtaining the corresponding central limit theorem by means of controlling the fluctuations of the endmost particles for obtaining there the exponential bound estimate sharpening of Lemma 7 in $\mathsection$\textit{4} below, and remark that the symmetry assumption is also employed through self-duality there. 


\vspace{3mm}  
\textit{\textbf{ 1.}} We set the notation and give necessary definitions to state the main results. The \textit{contact process} $\xi_{t}$ at rate $\mu$ on graph $G= (V,E)$ is a continuous-time Markov process on the set of subsets of $V$. Regarding sites in $\xi_{t}$ as occupied by, say, a particle, and other ones as empty, its transition rates may specified by the following set of rules: (i) a particle at $x$ gives birth to a new one at each empty $y$, $xy \in E$, at rate $\mu$, and (ii) particles die at rate 1. The process is said to be \textit{supercritical} when $\pr(\xi_{t} \not=\emptyset, \mbox{ for all } t)>0$, $\xi_{0}$ finite. For background on this extensively studied model, dating back to \cite{H74}, see \cite{D95} and \cite{L85, L99}. 

Let $Z_{M}$ be the graph with set of sites the integers $\Z$ for which different sites at Euclidean distance at most $M$ are adjacent\footnote{Our results and approach will be easily seen to apply under much less restrictive assumptions on the interaction, although only the induced class of processes considered is treated for notational convenience. Namely, they will be seen to apply to any self-dual process, that is, the case where rule (i) is such that births occur at rate which is any function of the distance among sites.}. Let $\xi_{t}^{0}$ be the supercritical contact process on $Z_{M}$ started from $\{0\}$, let also $r_{t} = \sup \xi_{t}^{0}$.
\begin{thm}\label{THEthm} 
\hspace{1mm} There exists strictly increasing $(\psi_{n})$ such that, on $\{\xi_{t}^{0} \not= \emptyset, \textup{ for all } t\}$, $(r_{\psi_{n}} - r_{\psi_{n-1}}, \psi_{n} - \psi_{n-1})_{n\geq 1}$ are i.i.d.. 
\end{thm}

We shall make extensive use of the celebrated Harris' graphical representation, cf.\ \cite{H78}, for constructing versions of the process from different starting sets on the same probability space. The reader is refer to, for instance,  $\mathsection$ 3 in \cite{D95}, pp. 126-128, for definitions and standard associated terminology. To state the next result, let $\xi^{\Zm}_{t}$ be the contact process on $Z_{M}$ started from $\Zm$, where $\Zm :=\{0,-1,\dots\}$, and let also $R_{t}=\sup \xi^{\Zm}_{t}$.

\begin{thm}\label{Rr}
$\pr(r_{t} = R_{t}, \textup{ for all } t)>0 $.
\end{thm}

The following is the new consequence of ergodicity of the process mentioned. 

\begin{thm}\label{gen}
Let $\xi_{t}^{O}$ be a supercritical, translation-invariant contact process started from a single site, $O \in V$. If, for any $F$ finite, $\xi_{t}^O \cap F  = \xi_{t}^{V} \cap F \textup{ for all large } t$, almost surely on $\{\xi_{t}^O \not= \emptyset \textup{ for all } t\}$, then
$\pr(\xi_{t}^{V} \cap F = \xi_{t}^{ F} \cap F \textup{ for all } t)>0$.
\end{thm}

The remainder is organized as follows. In $\mathsection$\textit{2} we prove Theorem 3, and state and prove Proposition 4 mentioned in the introduction. The proof of Theorem 1 is given in $\mathsection$\textit{4}.

\vspace{3mm}

\begin{proof}[\textbf{ 2. }Proof of Theorem \ref{gen}] From translation invariance and a standard restart coupling argument, see for instance the proof of Theorem 2.30 (a) in \cite{L99}, one has that the hypothesis of the theorem implies that, for any $F$ finite,
\begin{equation}\label{equi} 
\pr(\xi_{t}^F \cap F  = \xi_{t}^{V} \cap F \textup{ for all large } t | \mbox{ } \xi_{t}^{F} \not= \emptyset, \mbox{ for all } t) =1. 
\end{equation}

Let $B_{n}= \{ \xi^{V}_{s} \cap F =\xi^{F}_{s} \cap F, \mbox{ for all } s \geq n \}$, for integer $n\geq0$.  A realization of the representation is denoted by $\omega$. We write that for all $\omega \in E_{1}$, $\omega \in E_{2}$ a.e.\ to denote that $\pr(\{\omega: \omega \in E_{1}, \omega \not\in E_{2}\}) = 0$. From $(\ref{equi})$ we have that for all $\omega \in \{\xi_{t}^{F} \not= \emptyset, \mbox{ for all } t\}$ there is an $s_{0}$ such that $\omega \in \{\xi^{V}_{s} \cap F = \xi^{F}_{s} \cap F \mbox{, for all }s\geq s_{0}\}$ a.e..  Hence we have that $\displaystyle{\pr\left(\cup_{n\geq 0} B_{n}\right) = \pr(\xi_{t}^{F} \not= \emptyset, \mbox{ for all } t)>0}$ and, by contradiction, that there is $n_{0}$ for which $\pr( B_{n_{0}}) >0$.

Let $B_{n_{0}}'$ denote the event such that $\omega'\in B_{n_{0}}'$ if and only if there exists $\omega\in B_{n_{0}}$ such that  $\omega$ and $\omega'$ are identical realizations except perhaps from any $\delta$-symbols (death events) in $F \times (0,n_{0}]$. Further, let $D$ denote the event that no $\delta$-symbols exist in $F \times (0,n_{0}]$. 
By independence of the Poisson processes in the graphical representation and then because $B_{n_{0}}' \supseteq B_{n_{0}}$, we have that
\begin{eqnarray*}
\pr(B_{n_{0}}' \cap D) &=& \pr(B_{n_{0}}') \pr(D) \nonumber \\
&\geq& \pr(B_{n_{0}}) e^{-|F| n_{0}} >0,  
\end{eqnarray*}
where $|F|$ denotes the cardinality of $F$, because $B_{0} \supseteq B_{n_{0}}' \cap D$ the proof is completed from the last display. To prove that $B_{0} \supseteq B_{n_{0}}' \cap D$, note that if $\omega$ and $\omega'$ are identical except that $\omega'$ does not contain any $\delta$-symbols that possibly exist for $\omega$ on $F \times (0,n_{0}],$ then $\omega \in B_{n_{0}}$ implies that $\omega'\in B_{n_{0}}$ and indeed $\omega'\in B_{0}$. 
\end{proof}

\textbf{Remark.}  As noted in \cite{T11}, Remark 2,  an argument akin to that in the last paragraph above can be invoked to explicitly show the concluding claim in the proof of Theorem 3 in \cite{MS00}.  

\vspace{2mm}

Let $(\hat{\xi}_{t}^{A}; \hat{\xi}_{0}^{A} = A)$ be the supercritical contact process on $\Zmg_{M}$, the subgraph of $Z_{M}$ induced by the non-positive integers, $\Zm$, defined by the corresponding representation. 
\begin{prop}\label{cons}
$\pr(\hat{\xi}_{t}^{\Zm} \cap F = \hat{\xi}_{t}^{F} \cap F \textup{ for all } t)>0$, for $F\subset \Zm$ finite.
\end{prop}

\begin{proof}
Note that the assumption of translation invariance was not used in the last part of the proof of Theorem \ref{gen}, and hence that it suffices to verify (\ref{equi}) for the conclusion there to hold; we will show the next stronger statement by using the rescaling result of \cite{BG90} that is known to generalize easily in this case (equivalently, one may employ the result in \cite{DS87} equally well in this case and use that the process at criticality dies out, show in \cite{BG90}, in order to avoid strengthening the assumptions under which we work here). Let $1(\cdot)$ denote the indicator function.

\begin{lem}\label{SHAPE1}
There is an $a>0$ such that the set of sites $x$ such that $1(x \in \hat{\xi}_{t}^{F})  = 1(x \in \hat{\xi}^{\Zm}_{t})$ contains $[-a t,0] \cap \Zm$ \textup{for all large} $t$, almost surely on $\{\hat{\xi}_{t}^{F} \not= \emptyset, \textup{ for all } t\}$.
\end{lem}
\end{proof}

\begin{proof}[\textit{Proof of Lemma \ref{SHAPE1}}.]
We need only consider the case $F=\{0\}$. By additivity, as the arguments will then be seen to apply for any $F=\{x\}$, the extension to all finite $F$ follows easily. Let us write $\hat{\xi}_{t}(x)$ for $1(x \in \hat{\xi}_{t})$. By means of the rescaling results referred to above and following the arguments in $\mathsection$ 6 of \cite{DS87}, see also Proposition 3 in \cite{T15}, we have that there are $C, \gamma \in (0,\infty)$ such that, for any $x \geq -at$,
\begin{equation}\label{fds}
\pr(\hat{\xi}_{t}^{0}(x) \not= \hat{\xi}^{\Zm}_{t}(x), \hat{\xi}_{t}^{0} \not= \emptyset) \leq Ce^{-\gamma t}
\end{equation}
$t\geq0$. Note that the statement for integer times then follows from (\ref{fds}) and the 1st Borel-Cantelli lemma since $\sum\limits_{n\geq1} \pr\big(\bigcup_{x\geq -an} \hat{\xi}_{n}^{0}(x) \not= \hat{\xi}^{\Zm}_{n}(x) | \hspace{1mm}  \hat{\xi}_{t}^{0} \not= \emptyset, \mbox{ for all } t\big) < \infty$, where we first used that  $\pr(\hat{\xi}_{t}^{0} \not= \emptyset \cap  \{\hat{\xi}_{t}^{0} \not= \emptyset, \mbox{ for all } t\}^{c})$ is exponentially bounded in $t$, where the last standard result is proved by a well-known restart argument, see Theorem 2.30 (a) in \cite{L99}.

Let $B_t^{x}= \bigcup_{t\in (n,n+1]} \{\hat{\xi}_{t}^{0}(x) \not= \hat{\xi}^{\Zm}_{t}(x) , \hat{\xi}_{t}^{0} \not= \emptyset\}$, $x \in \Zm$. We have that
\begin{equation}\label{BteM}
\textstyle{ \pr(B_t^{x}) e^{-2M\mu - 2} \leq  \pr(\hat{\xi}_{n+1}^{0}(x) \not= \hat{\xi}^{\Zm}_{n+1}(x),  \hat{\xi}_{n+1}^{0} \not= \emptyset),}
\end{equation}
where this inequality follows from the strong Markov property by letting $t_{0}$ be the first time such that $B_{t}^{x}$ occurs and considering the event that: (i) no particles attempt to occupy $x$ during $[t_{0},t_{0}+1]$, (ii) the particle of $\hat{\xi}^{\Zm}_{t_{0}}$ at $x$ does not die until $t_{0}+1$ and, (iii) one particle of $\hat{\xi}_{t}^{0}$ does not die until $t_{0}+1$. From  (\ref{fds})  and (\ref{BteM}) the result follows as before by simply noting that $\left\{ \exists \hspace{0.3mm} t_{m} \uparrow \infty:\bigcup_{x\geq -at_{m}} \hat{\xi}_{t_{m}}^{0}(x) \not= \hat{\xi}^{\Zm}_{t_{m}}(x)\right\}$ can also be written as $\left\{\exists \hspace{0.3mm} n_{k}  \uparrow \infty: \bigcup_{t\in (n_{k},n_{k}+1]} \bigcup_{x\geq -at} \{\hat{\xi}_{t}^{0}(x) \not= \hat{\xi}^{\Zm}_{t}(x)\}\right\}$. 
\end{proof}

\vspace{3mm}

\begin{proof}[\textbf{ 3. }Proof of Theorem \ref{Rr}] 
Let $\Zmg_{M}$ be the (full) subgraph of $Z_{M}$ induced by $\Zm$, where $\Zm$ denotes the non-positive integers. The proof of the next well-known result follows easily as in $\mathsection$ 2, (b), in \cite{DG83} by means of the extension of the key renormalization result there done in \cite{DS87}. The necessity of symmetry and translation invariance for the transition rates of the process in the arguments below is reflected by the crucial role this result plays.  

\begin{prop}\label{muM}
If the process is supercritical on  $Z_{M}$, it is supercritical on $\Zmg_{M}$.
\end{prop}
Let $\hat{\xi}_{t}$ be the process confined to $\Zm$ defined by the representation simply by neglecting birth arrows leading to sites $y \notin \Zm$.
Proposition \ref{cons} applied for $F= \M$ gives that 
\begin{equation}\label{nnnkuc1} 
\pr(\hat{\xi}_{t}^{\Zm} \cap \M = \hat{\xi}_{t}^{\M} \cap \M, \mbox{\textup{ for all }}t)>0. 
\end{equation}

Let $r_{t}^{\M} =\sup \xi_{t}^{\M}$, $\M := \{0,-1,\dots, -M-1\}$. It suffices to show that 
\begin{equation}\label{rtMcal}
\pr(r_{t}^{\M} = R_{t} , \mbox{ for all } t)>0,
\end{equation}
by considering $\{\xi_{1}^{0} \supseteq \M\} \cap \{  \xi_{s}^{0} \cap \{0\} \not= \emptyset \mbox { and } R_{s} \leq 0, \mbox{ for all } s \in (0,1] \}$, monotonicity and the Markov property, this is easy to see.

Let $Z_{M}^{1^{+}}$ be the subgraph of $Z_{M}$ induced by $\Z^{1^{+}} := \{1,2,\dots\}$. Let $\tilde{\xi}_{t+1}^{1\times 1}, t\geq0,$ be the process confined to $Z_{M}^{1^{+}}$ started from $\{1\}$ at time $1$. Let $S$ be the event that: there is a birth arrow from $\M \times [0,1]$ to $\{1\} \times [0,1]$ and $\tilde{\xi}_{t+1}^{1\times 1} \not= \emptyset, \mbox{ for all } t\geq0,$ and no $\delta$-symbols (death marks) appear in $\{1\} \times [0,1]$.  Let further $D$ be the event that no $\delta$-symbols exist in $\M \times [0,1]$. Note that on $S \cap D$ we have that $\{r_{t}^{\M} \geq 0, \mbox{ for all } t\geq0\}$. Letting $C= \{\xi_{t}^{\M} \cap \M  \supseteq \xi_{t}^{\Zm \backslash \M} \cap \M, \mbox{\textup{ for all }} t\}$ and using additivity, we have that, 
\begin{equation}\label{Cr}
\{r_{t}^{\M} = R_{t} , \mbox{ for all } t\} \supseteq  C \cap S \cap D .
\end{equation}

Let $C' = \{\hat{\xi}_{t}^{\M} \cap \M  \supseteq \hat{\xi}_{t}^{\Zm \backslash \M} \cap \M, \mbox{\textup{ for all }} t\}$. First by coupling and then by monotonicity we have that 
\begin{equation}\label{Csupsets}
C=\{\xi_{t}^{\M} \cap \M  \supseteq \hat{\xi}_{t}^{\Zm \backslash \M} \cap \M, \mbox{\textup{ for all }} t \} \supseteq C'.
\end{equation}

From $(\ref{rtMcal})$, $(\ref{Cr})$ and $(\ref{Csupsets})$ give that showing that $\pr(C'\cap S \cap D)>0$ suffices for completing the proof.  By (\ref{nnnkuc1}), monotonicity and the Markov property give that $\pr(C' \cap D)>0$. Proposition \ref{muM}, by translation invariance and the Markov property, gives that $\pr(S)>0$. However, $C' \cap D$ and $S$ are independent as are the Poisson processes that these events are measurable with respect to, and the proof is complete. \end{proof}
 
\vspace{2mm}
\begin{proof}[\textbf{ 4.} Proof of Theorem \ref{THEthm}] Given a space-time point $x \times s$, let $\bar{\xi}^{x\times s}_{t+s}, t\geq 0,$ denote the process started from $\{y: y\leq x\}$ at time $s$ and let also $R_{t+s}^{x \times s}= \sup\bar{\xi}^{x\times s}_{t+s}$; furthermore let $\xi^{x\times s}_{t+s}, t\geq 0$, denote the process started from $\{x\}$ at time $s$, and let also $r_{t+s}^{x \times s} = \sup\xi^{x\times s}_{t+s}$. We write that $x \times s \mbox{ c.s.e.}$ for $R_{u}^{x \times s}= r_{u}^{x \times s}, \mbox{ for all } u\geq 0$, where the shorthand c.s.e.\ stands for "controls subsequent edges". 

Theorem $\ref{Rr}$ gives that $p:= \pr( 0 \times 0\mbox{ c.s.e.})>0$. From this and the next lemma the proof of Theorem 1 easily follows by letting $\psi_{n} = \inf\{t\geq 1+\psi_{n-1}: r_{t} \times t \mbox{ c.s.e.}\},$ $n\geq0$, $\psi_{-1}:=0$, and simple, well-known arguments (see Lemma 7 in \cite{MS00}).

\begin{lem}\label{restrt}
Consider the non stopping time $\psi = \inf\{t\geq 1: r_{t} \times t \mbox{ \textup{c.s.e.}}\}$. 
We have that $\psi$ and $r_{\psi}$ are a.s.\ finite conditional on either $\{\xi_{t}^{0}\not= \emptyset, \textup{ for all } t\}$ or $\{0 \times 0 \textup{ c.s.e.}\}$. 
\end{lem} 

\end{proof}

\begin{proof}[Proof of Lemma \ref{restrt}]
Define the processes $\xi_{t}^{n}, n\geq0,$ as follows. Consider $\xi_{t}^{0}:= \xi_{t}^{0 \times 0}$ and let $T_{0} = \inf\{ t: \xi_{t}^{0} = \emptyset\}$; inductively for all $n\geq0$, on $T_{n}<\infty$, let $\xi_{t}^{n+1}:= \xi_{t}^{0 \times T_{n}}, t\geq T_{n},$ and take $T_{n+1} = \inf\{t\geq T_{n}: \xi_{t}^{n+1} = \emptyset\}$. 

Let $r_{t}^{n} = \sup\xi_{t}^{n}$ and consider $r_{t}' := r_{t}^{n}$ for all $t\in [T_{n-1}, T_{n})$, where $T_{-1}:=0$. Let $\tau_{1}=1$ and inductively for all $n\geq1$, on $\tau_{n}<\infty$, let $\sigma_{n} := \sum_{k=1}^{n} \tau_{k}$ and $\tau_{n+1} = \inf\{t \geq 0: R_{t}^{r'_{\sigma_{n}} \times \sigma_{n}} > r_{t}^{r'_{\sigma_{n}} \times \sigma_{n}}\}$, while on $\tau_{n} = \infty$ let $\tau_{l} = \infty$ for all $l\geq n$. Let also $N = \inf\{n\geq1: \tau_{n+1} = \infty\}$. Since on $\{\xi_{t}^{0} \not= \emptyset, \mbox{ for all } t\}$, and on its subset $\{0 \times 0 \mbox{ c.s.e.}\}$, we have that $\psi = \sigma_{N}$ and $r'_{\sigma_{N}} = r_{\psi}$, it is sufficient to prove that $\sigma_{N}, r'_{\sigma_{N}}$ are a.s.\ finite.  

We prove the last claim. Note that, by translation invariance and independence of Poisson processes in disjoint parts of the graphical representation, we have that for all $n\geq1$ the event $\{\tau_{n+1}= \infty\}$ has probability $p$ and is independent of the graphical representation up to time $\sigma_{n}$. This and Bayes's sequential formula give that $\pr(N =n) = p(1-p)^{n-1}$ and, in particular, $N$ is a.s.\ finite. Thus also $\sigma_{N}$ is a.s.\ finite, which implies that $r'_{\sigma_{N}}$ is a.s.\ finite because $|r'_{t}|$ is bounded above in distribution by the number of events by time $t$ of a Poisson process at rate $M\mu$. This completes the proof. 
\end{proof} 
 
\vspace{4mm}

\textit{Acknowledgments:} The author is grateful to Tom Moutford for comments and for pointing out an error in chapter $\mathsection$ 4 of an early version of his thesis \cite{T11}.   Thanks also to Ronald Meester that first became conscious of the overviewed reference to Proposition 4.9 in \cite{T11}, in the proof  of $(4.2.1)$, also in \cite{T11}, and to Daniel Ahlberg for kind correspondence bringing this to his attention.

\end{document}